\documentclass{amsart}

\usepackage{amsmath, amsthm} 

\newtheorem{Theorem}[equation]{Theorem}
\newtheorem{theorem}[equation]{Theorem}
\newtheorem{Proposition}[equation]{Proposition}
\newtheorem{proposition}[equation]{Proposition}
\newtheorem{lemma}[equation]{Lemma}

\newtheorem{conjecture}[equation]{Conjecture}

\newtheorem{Remark}[equation]{Remark}

\font\Bbb=msbm10
\def\F{\hbox{\Bbb F}} 
 
 \numberwithin{equation}{section}

\let\define=\def

\define\({ \left( }
\define\){ \right) }
\define\[{ \left[ }
\define\]{ \right] }
\define\<{ \langle }
\define\>{ \rangle }
  
\DeclareMathOperator{\Ima}{Im}
\DeclareMathOperator{\Ker}{Ker}

\DeclareMathOperator{\AG}{AG}
\DeclareMathOperator{\GL}{GL}
\DeclareMathOperator{\AGL}{AGL}

\DeclareMathOperator{\Aut}{Aut}
\DeclareMathOperator{\Sym}{Sym}

\begin{document}

\title[Automorphism subgroups for designs with $\lambda=1$]
{Automorphism subgroups for designs with $\lambda=1$}
 
 \author{William  M. Kantor
}

 \address{U. of Oregon,
 Eugene, OR 97403
  \ and \ 
 Northeastern U., Boston, MA 02115}
 \email{kantor@uoregon.edu}
 
\begin{abstract}

Given  an integer $k\ge3$ and a group $G$ of odd order$,$  
if  there exists a $2$-$(v,k,1)$-design and  if $v$ is sufficiently large
 then there is such a design whose automorphism group has a
subgroup  isomorphic to $G$.  Weaker results are  obtained when $|G|$ is even.

 \end{abstract}
  
\maketitle 

\section{Introduction}
About 40 years ago 
Babai \cite[p.~8]{Ba} proposed the following
``subgroup problem":%
\smallskip
 
\begin{center}
\parbox{.87\linewidth}{PROBLEM 2.7. 
Prove for every $k \ge 3,$ that$,$ given a finite group $G,$
\vspace{2.5pt}
there is a BIBD of block size $k$ (a 2-$(v,k,1)$-design) $X$ such that
$G\,\widetilde \le \, \Aut X.$}
  \end{center}
  \smallskip
  
\noindent
 R. \!M. \!Wilson proved this when $k$ is a multiple of $|G|$\, \cite[p.~8]{Ba}; \
\cite[Theorem~12.1]{LW}
 contains this when  $k-1$ is a multiple of $|G|$. (These results are  also in \cite[p.~311]{Wi3}.)
 
In this note we will  prove other  special cases of Babai's problem:

\begin{theorem}
\label{Babai odd} 
Given  an integer $k\ge3$ and a group $G$ of odd order$,$ 
if $v$ satisfies the divisibility conditions  for a $2$-$(v,k,1)$-design  and is sufficiently large  
then there is  a $2$-$(v,k,1)$-design whose automorphism group has a
subgroup  isomorphic to $G$.
\end{theorem}

When $k=3$ stronger results appear in \cite{Ca} and \cite{DK}.

\begin{theorem}
\label {Babai even}
Given an  odd integer $k\ge3$ and a group $G$ of even order
such that  ${(k , |G|)=1,}$ there are infinitely many $v$  for which there is a
$2$-$(v ,k,1)$-design whose automorphism group has a subgroup  isomorphic to $G$.
 \end{theorem}
 
\begin{theorem}
\label {Babai even k}
Consider an even integer   $k \ge4 $    and    a group $G$   of even order.
Assume that every prime power dividing $|G|$
either divides $k$ or is relatively prime to ${k(k-1)}$.
Then there are infinitely many $v$  for which there is a $2$-$(v ,k,1)$-design 
whose automorphism group has a subgroup  isomorphic to $G$.
 \end{theorem}
 
When $k$ or $k-1$ is a prime power, see~\cite[p.~8]{Ba} or \cite{Ka} for a stronger type of result: there  are infinitely many  $2$-$(v,k,1)$-designs  ${D}$  for which $G\cong\Aut {D}$.
Babai \cite[Conjecture~2.8]{Ba} asked for 
such a stronger result for arbitrary $k\ge3$,  but this presently seems out of reach:
there appears to be no method for  recovering the  classical  geometry 
underlying one of our designs as was done when $k$ or $k-1$ is a prime power.
See  Remark~\ref{recover A} and Section~\ref{Conjectures}
for further comments about proving this stronger result.

The single idea behind the above three theorems is to 
place copies of  a  $2$-$(p,k,1)$-design  in the lines of an affine space 
$AG(d,p)$  in a $G$-invariant manner, for large $d$ and a suitable prime $p$; this 
occurred in  \cite[Sec.~III.C]{Ka1} for a very different purpose.   
The  $2$-$(p,k,1)$-designs we use admit   suitable automorphism groups (which are cyclic 
for Theorems~\ref{Babai even} and \ref{Babai even k}), and 
are   special cases of lovely  results in \cite{Wi,Wi3, LW}. 

Theorem~\ref{Babai even} is proved in Section~\ref{Even},
while Theorem~ \ref{Babai even k} is in Section~\ref{Nets and  even $k$}.
The remainder of this paper is devoted to Theorem~\ref{Babai odd}: Section~\ref{odd} 
contains a proof that there  are  infinitely many 
designs behaving as in Babai's problem  when $|G|$ is odd,  while
Propositions~\ref{induction}  and \ref{Large designs} (based on 
Theorem~\ref{Moore and Ray-Chaudhuri})
contain the background needed for the proof of Theorem~\ref{Babai odd} at the end of  Section~\ref{More oddness}.

All of  our  proofs are the same for  abelian  and nonabelian groups.
In all of the results mentioned above $|G|$ is tiny relative to $v$.
Our theorems  do not deal with the case  $|G|\equiv 0$ (mod~4)  and $k\equiv 2$ (mod~4).
The case  $ (|G|,k)\ne1\ne  (|G|,k-1) $ seems especially difficult when $|G|$ is even.
\medskip

{\noindent \em Preliminaries}:  If $G$ is a group of permutations $x\mapsto x^g$ 
of a set $X$, and $L\subseteq X$, then 
$G_L:={\{g\in G\mid L^g=L\}}$ is the set-stabilizer of $L$ in $G$, which induces  the subgroup
$G_L^L$   of the symmetric group $\Sym(L)$.

A permutation group $C$ on a set $X$ is  {\em semiregular}
 if $x^c\ne x$ whenever $x\in X$ and $1\ne c\in C$;
 and $C$ is {\em regular} if it is transitive and semiregular.
 If  $\< c\>$  and $\< c'\>$ are  semiregular
 cyclic groups of permutations of $X$ having the same order then  $c$ and $c'$ 
 are conjugate in   $\Sym(X)$.   
 
 We will use the same symbol to denote a design and its set of points.

\section{Odd order} 
\label {odd}
\begin{theorem}  [Wilson \mbox{\cite[pp.~22-26]{Wi}}]
\label{Wilson}
Given $k\ge3 ,$ for all sufficiently large  primes $p\equiv 1$ {\rm (mod~$k(k-1)$)} there is a 
$2$-$(p,k,1)$-design {E} whose set of points is $F:=\F_p$ and whose
automorphism group  contains  $\{x\mapsto  x+b\mid b\in F\}$.

Moreover$,$  if $p=1+k(k-1)t$  with  $t$  odd  then {E} can be chosen so that
$\{x\mapsto s x \mid s\in F, \,s^{t}=1\}$ is also a group of automorphisms of {E}.
\end{theorem}

If $t=(p-1)/k(k-1) $ is odd, the subgroup  of $F^*$ of order
 $2t$ factors as $S\times \<-1\>$ for a subgroup $S$ of order  $t$.
 Then \cite{Wi} obtains $A\subset F$  such that $\{ sA+b \mid s\in S, \, b\in F \}$ is the set 
of blocks  of $E$.  

 The preceding theorem lets us handle Babai's problem when $|G|$ is odd:
 
\begin{theorem}
\label{Babai odd prime power} 
Given  an integer $k\ge3$ and a group $G$ of odd order$,$ 
 there are infinitely many $v$  for which there is a $2$-$(v ,k,1)$-design 
whose automorphism group has a subgroup  isomorphic to $G$.
\end{theorem}

\proof
By Dirichlet's Theorem  there is a  prime $p\equiv 1 +{k(k-1)|G|}$ (mod~${2k(k-1)|G|)}$.
  If we write $p-1= k(k-1) t$,
 it follows that  $(p-1)/\{k(k-1)\}=t$  is odd  and divisible by  $|G|$. 
 As above, let $F=\F_p$ and let $S$ be the subgroup of $F^*$ of order~$t$.
  
We will prove the theorem by using suitable powers  $v=p^d$. 
Let  $V=F^d$, where $d$ is chosen so that $G$ is (isomorphic to) a group of permutations of a basis of $V$ and hence is in $\GL(V)$. (For example,  any integer $d\ge |G| $ can be chosen.)

We will use the affine space ${\bf A}:=\AG(d,p)$  whose set of points is  $V$.
Clearly $G<\GL(V)<\AGL(V)$. 
(Here $\AGL(V)=\{v\mapsto vM+c\mid M\in \GL(V), \, c\in V\}$ is 
$\Aut {\bf A}$ if $d>1$.)
Let $\mathcal L$ be a set of representatives of the orbits of $G$ 
on the lines of  ${\bf A}$.%

Let  $L\in \mathcal L$.  View $L$ as $F$, so the group $\AGL(1,p)$ of $p(p-1)$  affine transformations $x\mapsto ax+b$ for $a\in F^*$, $b\in F$, 
corresponds to the affine group $\AGL(L)$ on $L$  obtained from  $\AGL(V)$.
Then  $\{x\mapsto s x +b\mid s\in S, b\in F\}$ corresponds to a subgroup $S (L)$
 of  $\AGL(L)$ of order $p t$.
 Each subgroup of $\AGL(L)$ of order dividing $|S|=t$ lies in $S (L)$ 
 (since the quotient group $\AGL(1,p)/\hspace{-.5pt}\{x\mapsto  x+b\mid b\in F\}$
 is isomorphic to the cyclic group $F^*$).

The  set-stabilizer $G_L$ induces on $L$ a subgroup $G_L^L$ of $\AGL(L)$.
Since $|G|$ divides $t=|S|$ so does $|G_L^L|$.  Then $ G_L^L  \le S (L)$
 by the preceding paragraph.  (In fact, $G_L^L$ is even more restricted since $p>|G_L^L|$, but we will not need this fact.)
 
Use each $L\in \mathcal L$ as the set of points of a  
$2$-$(p,k,1)$-design $ {D}_L $ behaving as $E$ does at the end of Theorem~\ref{Wilson},
so $G_L^L  \le S(L) \le \Aut  {D}_L$.
(The end of  Theorem~\ref{Wilson} required that   that $t$ and $|G|$ are odd.)

For each $L\in \mathcal L$ let ${\mathcal B}_L$ be the set of blocks of~${D}_L $.
 If $g\in G$
let  $ {D}_{L^{{}^g}}$  denote the design  $ ({D}_L)^g$ whose set of points is  $L^g$
and whose set of blocks is $({\mathcal B}_L)^g$.

 {\em This well-defined}: {\em  if $L^g=L^{g'}$ for $g,g'\in G$ then 
$({D}_L)^g=({D}_L)^{g'}$.}  For, if $h=g'g^{-1} $ then 
$h \in G$ and $L^h=L$, so  the  permutation $h^L$ induced by $h$ on $L$ lies in  
 $  G_L^L\le \Aut  {D}_L$.  Then $h^L$  sends  $ {D}_L $ to 
itself,  so $({D}_L)^g=({D}_L)^{g'}$, as required.

Define a  design ${D}$ as follows:  

\hspace{20pt} points are the points of {\bf A} 

\hspace{20pt} blocks are the elements of 
$ \! \! \displaystyle\bigcup _{L\in \mathcal L, \, g\in G} \!  \!  \!  \!  \! ({\mathcal B}_L)^g$.
 \vspace{1pt}
 
\noindent
It is elementary that ${D}$  is  a 2-$(p^d,k,1)$-design:  any two points lie in a unique line
$L^g$ for $L\in {\mathcal L}$ and $g\in G$, and then in a unique member of 
$({\mathcal B}_L)^g$. Since  $G$ is in $\AGL(V)$ and permutes the sets $({D}_L)^g$ 
 it is a subgroup of  $\Aut {D}$. \qed
 
\vspace{1pt}
\Remark\rm
By the last sentence of Theorem~\ref{Wilson}, 
the first paragraph of the above proof contains a solution to
 Babai's problem for the cyclic group of order $|G|$. 
 The proofs of Theorems~\ref{Babai even}  and  \ref{Babai even k} involve  something similar:
 a cyclic group case of  Babai's problem is used to deal with much more general groups.
 
\Remark\rm
Placing designs on the blocks of another design is standard \cite[p.~28]{Wi}.  
Preserving the automorphism group is less standard.  
The above simple  method was used  in  \cite[Sec.~III.C]{Ka1}
to construct flag-transitive designs;  preserving a group of  automorphisms
of the larger design was as essential there as it  is here.

\Remark \rm
\label{recover A}
We used {\bf A}  with an arbitrary group  $G$ of odd order. 
Given the action of $G$ on $V$, the groups  $G_L$  and  $G_L^L$ are known;
since $p>|G|$, the group  $G_L^L$ is cyclic.

However, there is  flexibility with the designs ${D}_L$.  
We only needed to have $G_L^L\le \Aut {D}_L$ (for each $L\in  \mathcal L$)
in order for the proof to work.  Thus, each of the original designs 
${D}_L$ ($L\in  \mathcal L$) can be replaced by $({D}_L)^{h(L)}$ for any permutation
$h(L)$ of the points of $L$ that normalizes $G_L^L$.

Suitable changes of this sort  might provide a way to 
obtain a 2-$(p^d,k,1)$-design ${D}'$ such that $G\cong \Aut {D}'$. For this purpose 
{\em it appears to be necessary to  recover the affine space {\bf A}  
from  some such design ${D}'$}.  However, we have been unable to do this
(cf. Section~\ref{Conjectures}).
 
 \Remark \rm
 On the other hand,  {\em each design ${D}_L$ admits the group $S(L)<\AGL(L)
 =\AGL(V)_L^L$  as a  group of automorphisms that is {\em regular} on blocks}:  
 $\{{sA+b\mid}    s\in S, \, b\in  F\}$ is the set of  $|S|p=p(p-1)/k(k-1)$ blocks of each design 
 constructed in \cite[p.~22]{Wi} starting from a suitable initial block 
 $A\subset F$.
 Once again this uses the fact that $t$ and $|G|$ are odd.
  
   \Remark 
 \label{identity}
 If $B$ is a block of the design $D$  constructed in the  proof 
 of {\,\rm Theorem~\ref{Babai odd prime power}} then $G_B^B=1$.
 \rm  For, $B$ is in a  unique line $L$ of {\bf A}, so $L$ is fixed by $G_B$.
 Then $G_B^L\le S(L) $  as in the above proof. However, as already noted in the preceding remark, $S(L)$ is regular on the blocks of  ${D}_L$,  so $G_B^L \le S(L) _B =1 $ 
 and hence~$G_B^B =1 $.
   
 {\em This will be crucial in}  Section~\ref{More oddness}.
 
 Note that $G$ has many  fixed points,  so   there are
  many lines of {\bf A}  fixed pointwise by many elements of $G$.
  
\section{Theorem~\ref{Babai even}}
\label{Even}
When $|G|$ is even we use a consequence of a theorem 
of  Lamken and Wilson \cite{LW}; but first we need a prime:

\begin{lemma}
\label {Lamken Wilson lemma} 
Let $k\ge3,$ and let $h$ be a multiple of $4$  such that $(k, h )=1$.
 Then there are infinitely many primes $p>h$ satisfying the following conditions
 for some integer $n$$:$
\begin{itemize}
\item[(i)]$p=1+(k-1)n,$

\item[(ii)]$n(n-1)\equiv 0 $ {\rm(mod $k$),} 

\item[(iii)]$n(n-1)\equiv 0 $ {\rm(mod $4k$)}  if $k\equiv 3 $ {\rm(mod 4),} and 

\item[(iv)]$(p-1,h)=(k-1,h).$
\end{itemize}
\end{lemma}

\proof 
Let  $w$ be a positive integer such that $ k w\equiv 1$ (mod $h  $). Then
$\big(1+k(k-1)w, {  h k (k-1)\big)}=  
\big({1+k(k-1)w,} h \big)=(1+{(k-1),} h  )=1. $
By Dirichlet's Theorem there are infinitely many integers $y$ such that  $p:=1+k(k-1)w + \{    h k  (k-1)\}y =1+(k-1) n$ is prime, where  $n:= k w+   h ky \equiv 0$ (mod $k$).
Then (ii) is clear, and  (iii) holds:
  $n -1=(k w-1)+ h ky $  is a multiple of $h$ and hence of 4.
 Finally, (iv) holds:  $(p-1,h)=\big(k(k-1)w +  \{ h k (k-1)\}y,h\big)  
=\big(kw(k-1)  ,h\big) =(k-1  ,h) $.  \qed

\begin{Theorem}
[Lamken and Wilson \mbox{\cite[Theorem~12.1]{LW}}]
\label {Lamken Wilson} 
Given $k,$ for all sufficiently large $p$ satisfying the first three conditions
of \hspace{.7pt}{\rm Lemma~\ref{Lamken Wilson lemma}} there is a
$2$-$(p,k,1)$-design $E$ such that $\Aut E $ has a cyclic subgroup of order $k-1$ 
having one  fixed  point and semiregular on the  remaining points.

\end{Theorem}

{\noindent\em Proof of }Theorem~\ref {Babai even}.
We imitate the proof of  Theorem~\ref{Babai odd prime power}.  
In Lemma~\ref{Lamken Wilson lemma} let $h:=|G|$,
where we increase $G$  if necessary in order to have $h$  divisible by 4.
(Admittedly this is  annoying.)  
Choose a sufficiently large $p>|G|$ so that the lemma applies.
Choose  $d$ sufficiently large so that $G$ is (isomorphic to)
 a subgroup of the symmetric group $S_d$ and hence  also of $\AGL(d,p)$.  
The points of our design ${D}$ are the points of  ${\bf A}=\AG(d,p)$.  

Let $L\in \mathcal L$, where $\mathcal L$ is a set of 
representatives of the orbits of $G$ on  the lines of  ${\bf A}$. 
Then $G_L^L\le \AGL(L)\cong \AGL(1,p)$ and $p>|G|\ge |G_L^L|$, so
$G_L^L$ is a cyclic group of order dividing $(p-1,|G|)= (k-1,|G|)$ 
by Lemma~\ref{Lamken Wilson lemma}(iv).
This  cyclic group fixes a point, and all remaining orbits have length 
$|G_L^L|$;  all permutations of $L$ having this cycle 
structure are conjugate in  $\Sym(L)$.
After identifying $L$ with the set of points of the 
design in Theorem~\ref{Lamken Wilson}  and conjugating by an element of $\Sym(L)$, we may assume that  $G_L^L$ is contained in the cyclic group of order $k-1$  provided by Theorem~\ref{Lamken Wilson}.  Thus,  $L$  is the set of points of  a design
${D}_L$,  isomorphic to the design $ { E}$ in that theorem, such that
$G_L^L \le \Aut {D}_L$.   

Now repeat the last three paragraphs of the proof of Theorem~\ref{Babai odd prime power}.\qed
  
 \section{Moore and Ray-Chaudhuri}
 \label{Moore}
 Wilson \cite[p.~29]{Wi} credits Ray-Chaudhuri for the following
 generalization of a standard, fundamental  result due to Moore  \cite[p.~276]{Mo}:

  \begin{theorem}
\label{Moore and Ray-Chaudhuri}
A $2$-$(w,k,1)$-design $W,$
  a transversal design $TD(k,y-x)$  
  and a  $2$-$(y,k,1)$-design $Y$ with an $x$-point subdesign $X$ 
  produce a $2$-$(w(y-x)+x,k,1)$-design.
\end{theorem}

Here a  {\em transversal design} $TD(k,n)$ consists of $kn$ points, 
 $n^2$  subsets  of  size $k$  called ``blocks'', and a partition
of the points into $k$ ``groups'' of size $n$, such that each block meets each group
in a single point and any two points in different groups are in a unique block.

The following  proof is based on \cite[pp.~29-30]{Wi}, and is 
included since we need properties of the constructed design. 
\proof  
If $Z:=Y-X$ as a set of points,  then 
$X\cup (W\times Z)$ will be the set of   points of our  new design.
Let $A$ be a block of $W$, hence of size $k$.  There is a  transversal  design 
$TD(k,y-x)$  on $A\times Z$ whose set of groups is $\{ a\times Z\mid a\in A \}$
and whose set of blocks will be denoted ${\mathcal B}_{A\times Z}$;
this transversal  design, denoted ${T }_{A\times Z}$,
has nothing to do with the design on $Y$. 

Imitating Moore \cite[p.~276]{Mo} produces a new design  as follows:

\ \ points:  elements of $ X\cup (W\times Z ) $;

\ \ blocks are of four sorts:

\begin{itemize}
\item the blocks of $X$,
\item  for each $a\in W$ and each block $B$ of $Y$ not inside $X$, 
\begin{itemize}
\item
$a\times B$ if $B\cap X=\emptyset$, or 
\item  
$ x\cup \big (a\times (B- x)\big ) $ if $B\cap X =x $,
 and
\end{itemize}
\item  $\bigcup \{{\mathcal B}_{A\times Z} \mid \mbox{$A$  is a block  of $W$}\}$.
\end{itemize}

There is no conflict between the blocks  in ${T }_{A\times Z}$ for different choices of $A$: distinct intersecting sets $A\times Z$  and $A'\times Z$  intersect in 
a group $a\times Z$ that meets each block of ${T }_{A\times Z}$ or 
${T }_{A'\times Z}$   once.

 The only other part  the proof worth a comment concerns a
 pair  $(a_1,z_1),(a_2,z_2)\in W\times Z$ with $a_1\ne a_2$.  
 Since $a_1\ne a_2$ there is a unique block  $A$  of $W$ containing them, and  
  $(a_1,z_1) $ and  $(a_2,z_2)$ belong to different groups 
  $a_1\times Z$ and  $a_2\times Z$  of
  ${T }_{A\times Z}$.  Then  there is a unique block
in  ${\mathcal B}_{A\times Z} $  containing them.  \qed 

\Remark
\label{n_0}\rm
The existence of a  $TD(k, n) $ 
 is equivalent to the existence of a set of $k-2$ mutually orthogonal Latin
squares of order  $n$  \cite[Lemma~2.1]{Wi2}.  If $N(n)$ denotes the maximum number of 
mutually orthogonal Latin squares of order  $n$, then \cite{CES} proves that
there is an integer $n_0$ such that $\, N(n)\ge \frac{1}{3}n^{1/91}\, $
if $\, n>n_0$ (and there  are better bounds known \cite{Wi2}).
  Thus, if $n(k):= \max(n_0, (3k)^{91})$ then  
  \begin{equation} 
\label{TD}
\mbox{\hspace{-140pt}  \em If  $n>n(k)$ then there is a $TD(k,n)$.}
\end{equation}

 \section{Nets and  even $k$}
\label{Nets and  even $k$}
As in Sections~\ref{odd} and  \ref{Even}
the proof of Theorem~\ref{Babai even k} requires a suitable 
design  on a prime number of points.  
Whereas Theorem~\ref{Babai even} used a $2$-$(p,k,1)$-design having a cyclic automorphism group of order $k-1$ fixing one point and semiregular on the remaining points 
(Theorem~\ref{Lamken Wilson}),
this time we need  a $2$-$(p,k,1)$-design having a cyclic automorphism group of order $k$ fixing one point and semiregular on the remaining points 
(Theorem~\ref {prime}).
For this purpose we use Theorem~\ref{Moore and Ray-Chaudhuri} and transversal designs.  However, it will be easier to start with  nets.
 \subsection{Nets}
 \label{Nets}
 The dual of a transversal design $TD(k,n)$ is a $(k,n)$-{\em net}:  a set of $n^2$ points and 
 $kn$ subsets  of  size $n$  called ``lines''
such that  distinct lines meet at most once and the points are partitioned
 into $k$ ``parallel classes" each consisting of $n$ lines.  (Parallel classes correspond to groups.)
 The examples we need arise from unions of $k$ parallel classes of lines of a desarguesian affine plane $\AG(2,n)$; the translation group of the plane acts as a group of automorphisms of the net. Clearly these 
 examples exist whenever $n$ is a prime power and $k\le n+1$.

\begin{lemma}
\label {act on F}

Let  $q $   and $ m$ be powers $>1$ of a prime $p,$  
 and $E=\F_{q^{ m}}\supset F =\F_{q}$.
Let $\sigma\colon\! x\mapsto x^{q}$ and let $T\colon E\to F $ be the trace map.
If $a \in  E-\Ker T $ and $h\colon x\mapsto x\sigma+a,$
then $\<h\>$ has order $ pm$ and is semiregular on $E.$  
\end{lemma}
 
 \proof 
By induction, $h^i\colon\! x\mapsto x {\sigma^i}+\sum_{j=0}^{i-1}a {\sigma^j}$ 
for all $i\ge 1$, so  $h^m\colon \!x\mapsto x {\sigma^m}+ T(a) =x  + T(a)$ and
$h$ has order  $p m $.

 If $x\in E$ then $T(x) :=\sum_{j=0}^{ m-1}x {\sigma^j}$ and
 $T(x(\sigma -1))= T(x) \sigma -T(x)=0$,
so $\Ima(\sigma-1)\subseteq \Ker T$.
If  $i\ge1 $,     $d:=(m ,i)<m$  and  $x { \sigma^i}=x$,  then ${T(x)=(m /d)
\sum_{j=0}^{ d -1}x {\sigma^{j  }}} =0$  since $ m /d$ is a multiple of $p$,
so  $\Ima(\sigma-1)+\Ker(\sigma^i -1)\subseteq \Ker T$.

For semiregularity,   let $0<i<p m $ and suppose that $h^i$ fixes $x$.  Then 
$x{(1-\sigma^i)}=\sum_{j=0}^{i-1}a {\sigma^i}$, so
$x(1-\sigma^i)(\sigma -1)=a({\sigma^i}-1)$ and  $x (\sigma -1) + a
\in \Ker({\sigma^i}-1)$.  If $i\ne m$ then $(m ,i)<m$; we have seen that this implies that
 $ \Ker T$   contains $  \Ima(\sigma-1)+\Ker({\sigma^i}-1)$    and hence also $a$, 
 which is not the case. 
If $i=m$ then we obtain the same contradiction: $0=x(1 - \sigma^m )
=\sum_{j=0}^{m-1}a {\sigma^j} =T(a).$
 \qed  

\begin{lemma}
\label {product}
Let  $q$ and $ m$ be powers  $>1$ of a prime $p$. Let  
  $3\le k<q$. Then there is a $(k,q^{m})$-net
having a cyclic automorphism group of order $ m $
that is semiregular on both the points and lines and leaves invariant  each parallel class.
 \end{lemma}
 
\proof 
We use the notation of the preceding lemma.
Consider the affine plane $\AG(2, q^{m})$ defined using $E$.
Our net will consist of the points of this plane and~any union of $k$ parallel classes of lines
of the form $y=t x+b$ with $1\ne t \in F $ (so $t \sigma= t$). 

Let $g\colon \!(x,y)\mapsto (x^h,y^h)$.
By the preceding lemma,  $\<g \>$  has order $ p m$ and is semiregular on points.
Moreover, if  $i\ge1$  then $g ^i$  sends the line $\{(x,t x+b)\mid x\in F\}$ to the parallel line
$\{\big (x {\sigma^i} +a_i , t (x {\sigma^i})  +b {\sigma^i}+a_i\big)\mid x\in F\}$,
where $a_i : =\sum_{j=0}^{i-1}a {\sigma^j}$.    As above, 
$a _i (\sigma-1) = a(\sigma^i-1) $.

We still need semiregularity on lines. 
If   $0<i<p m$ and  $g^i$ fixes a line $y=t x+b$ of the net, then 
$t (x {\sigma^i})+b {\sigma^i}+a_i= t (x {\sigma^i}  +a_i) +b $, so 
$b(\sigma^i-1)=a_i(t -1)$.  Then $b(\sigma -1)(\sigma^i-1)=a(\sigma^i-1)(t -1)$, so  $b (\sigma -1)-a(t -1)\in \Ker(\sigma^i-1)$.
If $ i\ne m $ then 
 $a (t -1) \in \Ima( \sigma -1)+\Ker(\sigma^i-1)\subseteq \Ker T$  
(as seen above), 
which is impossible since $0\ne t -1\in F $ and $a$  is not in the  $F$-space $  \Ker T$.
Thus, $i=m$ and  $0=b(\sigma^m-1)=a_m(t -1)=T(a)(t -1)$, which is again impossible
since $t\ne1$.

This proves that  $\<g ^p\>$ behaves as required.
 \qed

\begin{lemma}
\label {s,t}
 Let $k\ge3$  be an integer    and let $p_1,\dots,p_r$ be its  distinct prime factors. 
 \vspace{2pt}
 For each $i$   let $m _i > k $  be a power of $p_i,$  
so  $k|\pi:=\prod_im_i$.

 Then for each  integer $ s > n(k)$ 
  there is a  $(k, s \prod_im_i^{m_i})$-net having a cyclic automorphism group
  of order $\pi  $
   that is semiregular on both points and lines while leaving invariant each parallel class.

 \end{lemma}
 
 Remark~\ref{n_0} contains the definition of $n(k)$.
 We emphasize that $s$ and the $m _i$ are not related. 

\smallskip
 \proof
For each $i$, by using Lemma~\ref{product} with $q=m=m_i$ we obtain a 
$(k, m_i^{m_i} )$-net $N_i$ having a cyclic automorphism group $C_i$ of order ${m_i } $ that is semiregular on both points and lines 
and leaves invariant each parallel class. 

 By  Remark~\ref{n_0}, if $s > n(k)$  then  there is a $(k, s )$-net $N_ \infty$.
The net required  in the lemma is a product  $N=N_1\cdots N_ r N_ \infty $, which we now define.

Let $X_i$ be the set of points of $N_i$ 
and  ${\mathcal L}_{i1},\dots,{\mathcal L}_{i k}$
 the parallel classes of $N_i$, 
 so $\bigcup_j{\mathcal L}_{ij}$ is the set of lines of $N_i$.  Then 
 $N$ is defined as follows:
$X:=X_1\times \cdots\times  X_r\times X_ \infty $ is its set of points, 
while its parallel classes are 
${\mathcal L}_j  :={\mathcal L}_{1 j}\times \cdots\times{ \mathcal L}_{r j}\times {\mathcal L}_{\infty j} $,
$1\le j\le k$, and $\bigcup_j{\mathcal L}_j  $ is its set of lines.  

In general, the groups $\Aut N_i$ are not involved in $\Aut N $ since we used  an arbitrary ordering  
of the parallel classes of each $N_i$.  However, for our purposes this is not a problem since $C_i$ leaves invariant each parallel class of $N_i$ and we will use the identity on $N_ \infty$:
there is a cyclic  automorphism group $C\cong \prod_iC_i  $ 
 \vspace{2pt}
 of  $N$ of order  $\pi=\prod_i   {m_i}  $, consisting of all  
 $(x_1,\dots,x_r,x_\infty)\mapsto (x_1^{c_1},\dots,x_ r ^{c_r} , x_\infty)$
for $c_i\in C_i$.    
(This is an automorphism of $N$:  $c_i$ permutes the lines in each 
${\mathcal L}_{i j}$, so if 
$(x_1,\dots,x_r,x_\infty)\in 
(L_{1j},\dots , L_{rj}, L_{\infty j})\in {\mathcal L}_j$ then
$(x_1^{c_1},\dots,x_r^{c_ r},x_\infty)\in 
(L_{1j}^{c_1},\dots , L_{rj}^{c_r}, L_{\infty j})\in {\mathcal L}_j$.) 

The only  way a point  (or line) of $N$ can be fixed by the preceding element of  $C$ is for a point (or line) of every component to be fixed, 
and then each $c_i=1$ by the semiregularity of each $C_i$.
\qed
 
 \subsection{Theorem~\ref{Babai even k}}
 Following the models in Sections~\ref{odd} and \ref{Even} we need a prime $p$ and
 a $2$-$(p ,k,1)$-design admitting a suitable automorphism group. 
 
\begin{proposition}
\label {prime}
Given integers $k\ge 3$ and  $h\ge1 $  such that $(k-1,h)=1,$  
there are  infinitely many  primes  $p$ such that 
$(p-1,h)$ divides some power of $k$ and there is a $2$-$(p ,k,1)$-design 
having a  cyclic automorphism group  of order $k$ fixing one point and 
semiregular on the remaining points. 
 \end{proposition}

\proof 
We will use a design in Theorem~\ref{Moore and Ray-Chaudhuri} 
whose set  of points is $U:=X\cup\big (W\times Z\big),$  
$ Z:=Y-X,$
where $X $ is  a subset  of size 1 of the design $Y$  \cite[Corollary~2C.1]{BS}.
For this we need three ingredients involving one choice of a suitable prime $q>h$,
a suitable choice in Lemma~\ref{s,t} of  the $m_i$  such that $\pi=\prod_im_i$  
is divisible by $k$,  and   
infinitely many~$s$:%
\begin{enumerate}
\item[(1)]  a $2$-$(qk ,k,1)$-design $W$ having 
a cyclic automorphism  group $C$ of  order  $k$  that is semiregular on points and whose 
$q$ point-orbits are blocks  of $W$  \cite[Theorem~1.2 and p.  308]{Wi3},
  
 \item[(2)] a transversal design $T=TD(k, (k-1)h  \pi s)$
having a cyclic automorphism group of order $k$ that is semiregular on both points and lines
($T$ exists for all sufficiently large $s$ by Lemma~\ref{s,t}),
and 

\item[(3)] a  $2$-$(y ,k,1)$-design $Y$  with $y:=1+  k(k-1)   (\pi/k) s $,
   and an arbitrary point  $X$  of $Y$ \
($Y$ exists for all sufficiently large $s$ \cite[Theorem~1.1]{Wi}).
(Note that we do not have any information concerning automorphisms of $Y$.)

 \vspace{-2pt}
\end{enumerate}
Moreover, we require that
\begin{enumerate}
 \vspace{-2pt}
\item[(4)]   $p=1+qk\cdot (k-1) \pi s=1+|W|(y-1)$ is   prime$,$ and 
\item[(5)]   $(p-1,h)$ divides some power of $k$
\rm(a condition  in the proposition).
\end{enumerate}

We will proceed in four steps.

{\smallskip\noindent (I) }  \emph{Number Theory$\,:$    $\pi,$
$p$ and $s$.} 
Write $h=h_0h'$, where $(k,h')=1$  and all  primes dividing $h_0$ also divide $k$.
Then $(k(k-1),h')=1$ since $(k-1,h)=1$, so that $h'$ is odd.
 We may assume that  $h_0$ divides the product $\pi $  of  the   
 $m_i $ used in Lemma~\ref{s,t}.%
 
We have  a prime $q>h$ in (1),  so $(q k (k-1)\pi ,h')=1$  since 
$(k ,h ')=1=(k-1,h)$. 
Let $t$ be a  positive integer such that ${q k (k-1)\pi t}\equiv 1$ (mod~$h'$).
Then $(t,h')=1$   and
$\big (1+q k (k-1)\pi t,q  k (k-1)\pi h'\big)=\big (1+q k (k-1) \pi t, h'\big) =(1+1,h')=1 $
since $h'$ is odd.  
By Dirichlet's Theorem  there are infinitely many integers $f$ such that
 $$  p:=1+q k (k-1)\pi t+q k  (k-1)\pi h' f$$ is a prime.
 
Choose $s:=t+h' f$  with $f$ so large  that  the designs in (2) and~(3) exist.    Clearly (4) holds.

Moreover,   $(p-1,h)$ divides 
$\big(  q ({k-1})({t+h' f}) ,h'\big) k   \pi h_0 = \big((k-1) t,h'\big)k  \pi h_0= k  \pi h_0 ,$
 which divides some power~of~$k$,
 as required in (5).
 \smallskip
 
 Now that we have $W,$ 
 $T$ and $Y$  we need to turn the set  $U=X\cup\big (W\times Z\big)$ of size $p$ into a design.

{\smallskip\noindent (II) } \emph{The  cyclic group $\bar C$.}
 We need  \emph{a group $\bar C\cong  C$ of permutations of}  
 $U $.
Extend each  $c\in C$ (cf.~(1))  to  a permutation  $\bar c$ of $U$ that fixes the point in  $X$ and 
sends $(a,z)\mapsto (a^c,z)$ for $a\in W, \, z\in Z$.  Then $\bar C=\{\bar c\mid c\in C \}$ is~a group of  $k$
permutations of $U$ \emph{fixing  $X$ and semiregular on the remaining points.}  This is not yet a group of automorphisms of anything.

We will construct  the design in Theorem~\ref{Moore and Ray-Chaudhuri} by
 placing (in (III)) copies of the transversal design  
 $T$ in  the sets $B\times Z$ of size $k(y-1)$ arising from blocks $B$
of $W$, and  (in (IV))   copies of     
$Y$ in the sets $Y_a:=X\cup\big (a\times Z\big)$ of size $y$ 
for $a\in W$.

{\smallskip\noindent (III) }  \emph{Copies of   $ T$.}
We will use copies of the transversal design $T$  in (2) as  the transversal designs occurring  in the proof of Theorem~\ref{Moore and Ray-Chaudhuri}. 
(Recall that $y-1=k(k-1)h\pi s$.)
 
In view of the point-orbits in  (1), the stabilizer in $C$ of a block   of $W$ is either 1 or $C$.

Let ${\mathcal B}$ be a set of orbit representatives of $C$ on the blocks of $W$.
If $B\in {\mathcal B}$ let $T_{B\times Z}\cong T$ have $B\times Z$ as its set of points
and $\{b\times Z\mid b\in B\}$ as its set of groups.
If the stabilizer of $B$ in $C$ is 1,  $T_{B\times Z}$ is placed in $B\times Z$ arbitrarily.
 If $B$ is a $C$-orbit we have to be more careful.  Initially,   place
 $T_{B\times Z}$  arbitrarily.  We then have 
 two semiregular cyclic permutation groups  of order~$k$ on  $B\times Z$:
one is the restriction $\bar C^{B\times Z}$ of $\,\bar C$ to $B\times Z$, 
and the other is the cyclic automorphism group of $T_{B\times Z}$ provided by (2).  
These cyclic groups of order $k$ are conjugate 
by an element of  $\Sym\!\big (B\times Z \big)$;
conjugate by such an element in order to assume that $T_{B\times Z}$ has been placed in $B\times Z$ so  that   the 
cyclic groups coincide, and hence so that  
  $\bar C^{B\times Z}\le \Aut T_{B\times Z}.$

For $B\in {\mathcal B}$  and $c\in C$  let  $T_{{B^c\times Z}}$ denote  
the transversal design $  (T_{B\times Z})^{\bar c}$ having ${(B \times Z)^{\bar c}}=B^c\times Z$ as its set of points.  
This is well-defined:  if $B^c  =  B^{c'} $  
then $ {{\bar c' \bar c{}^{-1}} } $    induces  
the permutation $({ \bar c' \bar c{}^{-1}} ) ^{B\times Z}$   of $B\times Z$, which is
an automorphism  of $T_{B\times Z}$ by the preceding paragraph, so that 
$(T_{B\times Z})^{\bar c}= (T_{B\times Z}) ^{\bar {c}'}$.  

{\smallskip\noindent (IV) } \emph{Copies of  $Y$.}
 Next we place \emph{copies of the design   
$Y$ into the sets $Y_a =X\cup\big (a\times Z\big),$  $ a\in W,$} in the same manner.  Namely, let ${\mathcal W}$ be a set of orbit representatives of $C$ on the points  of $W$.
For $a\in {\mathcal W}$ place a copy of the design   $Y$ in $Y_a$ using
the bijection $X\mapsto X$,  $z\mapsto (a,z)$ with $z\in Z$;   then let
$Y_{a ^c} :=  (Y_a)^{\bar c}$ for $c\in C$.  As usual, this is well-defined since 
$a ^c=a^{c'}$ implies that  $c=c'$ by semiregularity (cf.~(2)).

 Using the construction in the proof of 
Theorem~\ref{Moore and Ray-Chaudhuri} we obtain a 
$2$-$(p,k,1)$-design   $U$ having $\bar C$ as a group of  $k$  automorphisms
that fixes the point $X$ and is semiregular on the remaining points. \qed

\medskip  

{\noindent \em Proof of} Theorem~\ref{Babai even k}. 
 Let $h=|G|$ and $p>h$ be as in the preceding proposition.
We imitate the proof of  Theorem~\ref{Babai odd prime power},
regarding  $G$ as a group of automorphisms of ${\bf A}=\AG(d,p)$
for   any sufficiently large $d$.  
Let $L\in \mathcal L$, where $\mathcal L$ is a set of~representatives of the orbits of $G$ on the lines of  ${\bf A}$.
Then $G_L^L$ has order dividing $\big ({p(p-1)}, |G| \big)=(p-1, |G|);$  by the proposition,  
this divides some power of $k$,  and hence divides $k$ by an hypothesis of the theorem.

The  cyclic group $G_L^L$ fixes a point  and is semiregular on the remaining points.
After identifying $L$ with the set of points of the 
design  $U$  in Proposition~\ref {prime} and conjugating by an element of $\Sym(L)$, we may assume that  $G_L^L$ is contained in the cyclic group of order $k$  provided by  the  proposition.  Thus,  $L$  is the set of points of  a design
${D}_L$,  isomorphic to   $ {U}$, such that
$G_L^L \le \Aut {D}_L$.   

Now  complete the proof by repeating the last three paragraphs of the proof of Theorem~\ref{Babai odd prime power}.\qed  
 
\section{Large designs}
\label{large designs}
The Doyen-Wilson Theorem \cite{DW} states that, whenever $y\ge2x+1$
and there are Steiner triple systems on $y$ and $x$ points, there is a Steiner triple system on $y$ points having a subsystem on $x$
points.  The following is a significant generalization of that result
\cite{DLL}: 
  
\begin{theorem}
  \label{Dukes}
 If $k\ge3 $  then there is an integer $ x_0(k) >k$  such that$,$ if  $x>x_0(k),$  
$y > xk,$ \ $x-1\equiv y-1\equiv0$ {\rm (mod $k-1$)}  and
$x(x-1)\equiv y(y-1)\equiv0$ {\rm (mod $k(k-1)$),}
then there is a $2$-$(y,k,1) $-design having an $x$-point subdesign.
 \end{theorem}
  
We use this for a  result  concerning large designs ($n(k)$ appears in Remark~\ref{n_0}):

\begin{Proposition}
  \label{induction}
Let ${\mathcal S}$ be a  set of $2$-$(u,k,1)$-designs$,$ 
 let  $\bar{\mathcal S}$ be the set of all such $u$ that occur for ${\mathcal S},$
and let $w\in \bar{\mathcal S}$.
Assume that$,$  if $x$ and $y$ are as in {\rm Theorem~\ref{Dukes}} with $y-x>n(k),$
then $x+w(y-x)\in \bar{\mathcal S}$.

Then $\bar{\mathcal S}$ contains all sufficiently large $u$ satisfying the divisibility conditions for a $2$-$(u,k,1)$-design.
  \end{Proposition}
  
Our  argument imitates  \cite{Ca}.
Note that the hypothesis involves only 
the initial existence of one $w\in \bar{\mathcal S}$.
 
  \proof  For $x_0(k)$ in  Theorem~\ref{Dukes},
let $x_1 > x_0(k)>k$ be any representative for a congruence class
  (mod~${k(k-1)} $)  of integers such that there exists a  $2$-$(x_1,k,1)$-design.
  Consider any integer $a\ge w x_1n(k)\ge w x_1 $.
Choose $y-x=x_1k (k-1)a>a>n(k)$, so $y-x> n(k)$, and then choose $x=x_1+k(k-1)t $ with $0\le t<a$, so $x\ge  x_1 > x_0(k)$.
 Then $x$ and $y=x_1+{x_1k (k-1)a}+k(k-1)t$ satisfy $y>kx$.
 (For, since $t<a$ and $x_1>k$, we have
  $y-kx = (k-1)\big(-x_1+x_1ka-k(k-1)t   \big), $
 where
 $x_1(ka-1) -k(k-1)t >k(ka-1)-k(k-1)a >0$.)
 
 Since $x$ and $y$  satisfy  the divisibility
 conditions and the requirements $x>x_0(k)$ and $y>kx$  in Theorem~\ref{Dukes}, there is
  a $2$-$(y,k,1) $-design having an $x$-point subdesign.
Theorem~\ref{Moore and Ray-Chaudhuri} also needs a $TD(k,y-x)$, 
which exists since  $y-x> n(k)$  (cf.~(\ref{TD})).
By hypothesis,  Theorem~\ref{Moore and Ray-Chaudhuri} 
produces a $2$-$(u,k,1) $-design such that 
 \begin{equation}
\label{u}
\mbox{\hspace{-13pt} $u:=x+w(y-x)\in  \bar{\mathcal S}$,\,   with\,
 $u=x_1+wx_1k (k-1)a + k(k-1)t$.}
\end{equation}
 Here  $u-1 \equiv0$ {\rm (mod $k-1$)}  and
$u(u-1)\equiv  0$ {\rm (mod $k(k-1)$).} We  will  show that 
\emph{the set of all $u$ obtained in $(\ref{u})$ contains the set of 
all sufficiently large $u\equiv  x_1$ {\rm (mod~$k(k-1))$} satisfying these divisibility conditions.}
 
 Given $a$, we have  $y-x=x_1k(k-1)a $ and $x=x_1+k(k-1)t$.
  By choosing $t=0,\dots,a-1$, we realize
 $$
 u=x_1+wx_1k (k-1)a,  \dots, x_1+wx_1k (k-1)a+k(k-1)(a-1).
 $$
For $y-x=x_1k  (k-1)(a+1),$ we realize
 $$
 u=x_1+wx_1k (k-1)(a+1), \dots, x_1+wx_1k (k-1)(a+1)+k(k-1)a.
 $$
 In order not to leave any gaps, we require that these intervals abut or overlap. This occurs as long as
$ x_1+wx_1k (k-1)a+k(k-1)a\ge  x_1+wx_1k (k-1)(a+1),
 $
 that is, $a\ge wx_1$, which is a  condition  already satisfied by $a$. So we can achieve all sufficiently large $x\equiv x_1$ (mod~$k(k-1)$). 

 Now let $x_1>x_0(k)$ run through  a set of representatives for  the congruence classes
mod~${k(k-1)} $ that satisfy the divisibility conditions for a $2$-$(x_1,k,1)$-design.~\qed
 
 \section{Theorem~\ref{Babai odd}}
 \label{More oddness} 
  
We call an automorphism group of a design
1-{\em blocked}  if the set-stabilizer of any block is the identity on the block;
our basic example was in Remark~\ref{identity}.
This notion is preserved by the construction  in Section~\ref{Moore}:

  \begin{proposition}
  \label{Large designs}
  \vspace{1pt}
  Let $k\ge3$ and let $G$ be a $1$-blocked automorphism group
  of a $2$-$(w,k,1)$-design $W $.  Then 
  a $2$-$(y,k,1)$-design $Y$ with a subdesign $X$ on $x$ points$,$  
together with  a transversal design $TD(k,y-x),$  produce a $\,2$-$(w(y-x)+x,k,1)$-design
$U$ such that $G$ is isomorphic to a $1$-blocked subgroup of  $\Aut\hspace{.5pt} U$.
\end{proposition}
 
\proof 
 We use the construction and notation in the proof of  
 Theorem~\ref{Moore and Ray-Chaudhuri}.
 Each $g\in G$ induces~on $U$  the permutation $\bar g$ sending
 $b\mapsto b$ and  $(a,z)\mapsto(a^g,z)$ for $b\in X,$ $ a\in W, $
 $z\in Z$.  Clearly $ G \cong  \bar G:=\{\bar g \mid g\in G\}$.
 
 For each  orbit-representative $A$  of $  G$ on the blocks of $W$
 we have a transversal design ${TD}_{A\times Z}$ whose set of points is $A\times Z$
and whose set of groups is $\{ a\times Z\mid a\in A \}$.
If $ g\in  G$ then $(A\times Z)^{\bar g}=A^{ g}\times Z $ for a block $A^g$ of $W$;
  let ${TD}_{ A^g\times Z } :=({TD}_{A\times Z})^{\bar g}$.  As in 
the proof of Theorem~\ref{Babai odd prime power}  this is well-defined:  
if $A^{  g} \times Z= A^{g'}\times Z $ with $g,g'\in G$ then $A^g=A^{g'}$,~so $g'g{}^{-1}=1$
 on $A$ since $G$ is $1$-blocked,~and hence 
 $({TD}_{A\times Z})^{\bar g}=({TD}_{A\times Z})^{\bar {g}'} $
since $\bar g'\bar g{}^{-1}=1 $~on~$A\times Z$.

By the construction in Section~\ref{Moore}, 
each $\bar g$ permutes the designs ${TD}_{A'\times Z}$  with $A'$ a block of $W$, and 
is the identity on any other block of $U$ (i.\,e.,  a block of $X$, or else
$a\times B$  or  $b\cup \big(a\times (B-b) \big)$ if $a\in W$ and if $B\cap X=b$).  Thus, $\bar G\le \Aut\hspace{.5pt} U$.

We need to verify that $\bar G$ is 1-blocked.
Consider a block  $E$  of $U$  fixed by $\bar g\in \bar G$.  By 
Section~\ref{Moore}, either $E$  is contained in  $X\cup(a\times B)$ for $a\in W$ and a block $B$ of $Y$, or $E$ is a block  of some ${TD}_{A\times Z}$.
In  the former case it is clear that $\bar g=1$ on $E$, 
 so we are left with $E$ in ${TD}_{A\times Z}$.  
 In view of the construction in  Section~\ref{Moore},
 $A$ is uniquely determined by $E$ and hence is fixed by $g$.
Since $G$ is $1$-blocked on $W$, it follows that  $g= 1$ on $A$.
Then $\bar g=1$ on $A\times Z$ and hence on $E$.
Thus, $\bar G$ is  a $1$-blocked subgroup of  $\Aut\hspace{.5pt} U$.  \qed
  
\Remark
\label{not even}
 An automorphism group of even order cannot be $1$-blocked. \rm
For,  an involution interchanges two points, hence fixes the block containing them 
and acts nontrivially on that block.
  \smallskip \smallskip
  
{\noindent\em Proof of} Theorem~\ref{Babai odd}. 
Apply Proposition~\ref{induction} to  the set  $\mathcal S$  of  $2$-$(v,k,1)$-designs 
 whose automorphism group has a $1$-blocked subgroup  isomorphic to $G$.
By Theorem~\ref{Babai odd prime power} and 
Remark~\ref{identity}, ${\mathcal S}$ contains some $2$-$(v,k,1)$-design.  
 
 We defined $n(k)$ in Remark~\ref{n_0} and  $x_0(k)$  in Theorem~\ref{Dukes}.
 Let $x>\max(n(k),x_0(k))$  and $y>kx $
 be integers such that there are $2$-$(x,k,1) $- and $2$-$(y,k,1) $-designs.  
 By  Theorem~\ref{Dukes}  there is a  $2$-$(y,k,1) $-design having an $x$-point subdesign.
 Since  $y-x>kx-x>n(k) $  there is  a  $TD(k,y-x)$ by (\ref{TD}).
Then  $x+v(y-x)\in \bar{\mathcal S}$ by Proposition~\ref{Large designs}.
 \
 Now use  Proposition~\ref{induction}.\qed
  
 \section{Conjectures}
  \label{Conjectures}
 Our theorems are significantly weaker than the corresponding results in 
 \cite{Ba, Ka,DK}, where $G$ is isomorphic to the {\em full} automorphism group of the constructed design.  We conclude with a conjecture concerning affine spaces that would produce designs with this stronger property.
 
\conjecture
\label{conjecture}
{\noindent \rm Given} an integer $  s\ge 14,$   a prime 
$ p\equiv 1\, ($mod~$s), $ and an affine space ${\bf A}'$ having the same set of points as the 
original affine space ${\bf A}=\AG(d,p), $ 
such that 

\qquad\qquad
for any subspaces $X $ of ${\bf A}$ and $Y'$ of ${\bf A}', $ 

\qquad\qquad
 either 
$X\cap Y'=\emptyset$ ~or~$|X\cap Y'|\equiv 1$~$($mod~$s)$.\medskip

{\noindent \rm Conjecture:}  
${\bf A}={\bf A}'.$ 
 
\smallskip
\rm

Note that it is essential here that $p$ is prime.  For suppose that $p=p_0^e>p_0$ 
for a prime $p_0\equiv 1$~$($mod~$s)$.  Let ${\bf A}_0 =\AG(e d,p_0)$,
 let {\bf A} be the set of affine $\F_p$-subspaces of ${\bf A}_0$.  
 If $g\in \AGL(d n,p_0)- \AGL(d,p)$ then $ {\bf A} ' :={\bf A} ^g $ 
 provides a counterexample to the conclusion  in the preceding conjecture.
 
 The condition $s\ge 14$ reflects the fact that 14 is the smallest integer $s=k-1\ge2$ 
 such that neither $s$ nor $s+1$ is a prime power:  when one of these is a prime power the 
 desired result is already known  \cite{Ba, Ka}.

\begin{theorem}
Assume that the preceding conjecture is correct.
Under the hypotheses in any of 
{\,\rm Theorems \ref{Babai odd}--\ref {Babai even k}}$,$ for infinitely many $v$ there is a $2$-$(v,k,1)$-design $D$ such that $G\cong \Aut D$.
\end{theorem}

\proof
By \cite{Ba,Ka} we may assume that neither $k$ nor $k-1$ is a prime power,
so that $p>s:=k-1\ge 14.$
Each theorem in Section~1 uses $2$-$(p,k,1)$-designs constructed in 
Theorems~\ref{Wilson} or
 \ref{Lamken Wilson}, or Proposition~\ref {prime}.
In view of those constructions,  in the situation of any of the theorems
in Section~1, 
there are  $2$-$(p,k,1)$-designs $E_1,E_2,E_3$ with $\AG(1,p)$  as their set of points
such that there is no isomorphism between any two of these designs
that lies in  $\AGL(1,p)$.  (Namely, start with a design $E_1$, and   for
 $i=2,3$   apply an $i$-cycle  of points to the blocks of $E_1$ in order to obtain the blocks of~$E_i$.)%

Let $d > 4$ be as  in the proofs, so we are using ${\bf A}=\AG(d,p)$
 based on a $d$-space $V$.
Let $\{v_1,\dots,v_ d\}$ be a basis of $V$.  There is a  connected
 graph $\Gamma$ with vertex set 
$\{v_1,\dots,v_ d\}$ such that $G\cong \Aut \Gamma$.

Let $c:=\sum_ 1^d v_ t $.
Place  $E_1$ in each affine 1-space $\<v_i\>$, place 
$E_2$ in each affine 1-space $\<v_i+  v_j\>+c$ such that $\{v_i,v_j\}$ is an edge of 
$\Gamma$, and place $E_3$ in every other affine 1-space of ${\bf A}$.  
(Note that,   since $d>4$, if $\<v_i+  v_j\>+c$=$\<v_{i'}+  v_{j'}\>+c$
then $\{i,j\}=\{i' ,j' \}.$)

This produces a 
$2$-$(p^d ,k,1)$-design $D$ with $G$ (isomorphic to) a subgroup of $\Aut D$.
(Compare \cite{Ka}.)
 
Let $h\in \Aut D$ and  consider the affine space ${\bf A}'={\bf A}^h$.  
If $X$ and $Y'$ are subspaces of {\bf A} and ${\bf A}',$  respectively, 
and if $X\cap Y' \ne \emptyset,$  then $X\cap Y'$ is the intersection of subdesigns and so is a subdesign.  Then $|  X\cap Y' |\equiv 1$ (mod $s$).
Thus,  ${\bf A}^h={\bf A}$  by our hypothesis concerning Conjecture~\ref{conjecture},
so  $h$ is an automorphism of ${\bf A}$,
and hence permutes the lines of ${\bf A}$.   Then $h$ also permutes the designs we have placed inside these lines,  so $h$ permutes the 
lines $\<v_i\>$.  The  intersection of these lines is 0, so
$h$ is a linear transformation.
 By construction,  $h$ also permutes   the lines $\<v_i+ v_j\> +c$, 
 so it induces an automorphism of 
$\Gamma$, and hence agrees with some
$g\in  G=\Aut \Gamma$ in its action on  the lines $\<v_i\>$.  Now $h':=hg^{-1}$
fixes each line $\<v_i\>$, and hence is a diagonal transformation:  $v_t^{h'}=a_tv_t$
for some $a_i\in K^*$ and all $t$.  Since $h'$ fixes each vertex of $\Gamma$ it fixes each edge:
$\<v_i+ v_j\>+c=(\<v_i+ v_j\>+c)^{h'}=\<a_iv_i+ a_jv_j\>+ \sum_ 1^d a_ t v_ t $. 
It follows that all $a_ t =1$ for $t\ne  i,j$, so $h'=1$ since  $\Gamma$ is  connected, and then~$h\in G$.\qed 

\Remark\rm
The fact that $X\cap Y'$ is a subdesign imposes  arithmetic and structural conditions that
can be included in the hypotheses of the above  conjecture.

 \smallskip\smallskip 
{\noindent\bf Acknowledgements.} \rm
I am grateful to Peter Dukes for assistance with \cite{DLL},
and to Jean Doyen for many helpful comments.
This research was supported in part by a grant from the Simons Foundation.

\end{document}